\documentstyle[12pt,amssymb,graphicx]{article}
\begin{document}
\renewcommand{\thefootnote}{}
\date{}

\def\thebibliography#1{\noindent{\normalsize\bf References}
 \list{{\bf
 \arabic{enumi}}.}{\settowidth\labelwidth{[#1]}\leftmargin\labelwidth
 \advance\leftmargin\labelsep
 \usecounter{enumi}}
 \def\newblock{\hskip .11em plus .33em minus .07em}
 \sloppy\clubpenalty4000\widowpenalty4000
 \sfcode`\.=1000\relax}

\baselineskip=12pt
\title{\vspace*{-1cm}
{\large BRUNNIAN LOCAL MOVES OF KNOTS AND VASSILIEV INVARIANTS}
}

\author{
{\normalsize AKIRA YASUHARA }\\
{\small Department of Mathematics, Tokyo Gakugei University}\\[-1mm]
{\small Nukuikita 4-1-1, Koganei, Tokyo 184-8501, Japan}\\
{\small e-mail: yasuhara@u-gakugei.ac.jp}
}

\maketitle

\vspace*{-5mm}  
\baselineskip=15pt

{\small 
\begin{quote}
\begin{center}A{\sc bstract}\end{center}
\hspace*{1em} 
K. Habiro gave a neccesary and sufficient 
condition for knots to have the same Vassiliev invariants 
in terms of $C_k$-move. 
In this paper we give another geometric condition 
in terms of Brunnian local move. 
The proof is simple and self-contained.
%
%
%K. Habiro defined a local move, {$C_k$-move}, 
%of oriented links and by using his \lq clasper theory',  
%he proved that two knots are same up to $C_k$-moves if and only if
%they have the same Vassiliev invariants of order $\leq k-1$. 
%K. Taniyama and the author gave a diagram-oriented proof of 
%this theorem by translating \lq clasper' into \lq band description' 
%of knots. The author think the proofs are elementaly, but little complicated. 
%In fact, we have to prove several sublemmas and lemmas
%for the proof. Here by using a local move, $(k+1)$-component 
%Brunnian local move, instead of $C_k$-move in their arguments, 
%we show that knots are same up to $(k+1)$-component 
%Brunnian local moves if and only if
%they have the same Vassiliev invariants of order $\leq k-1$. 
%This also gives a relation between local moves and Vassiliev invariants,  
%and the proof is much simpler than that of the theorem for $C_k$-moves.
\end{quote}}


\baselineskip=16pt

\bigskip\noindent
{\bf 1. Introduction}

\bigskip
A {\it tangle} $T$ is a disjoint union of properly embedded 
arcs in the unit $3$-ball $B^{3}$. 
Here $T$ contains no closed arcs. A tangle $T$ is {\it trivial} if there 
exists a properly embedded disk in $B^3$ that contains $T$. 
A {\it local move} is a pair of trivial tangles 
$(T_{1},T_{2})$ with $\partial T_{1}=\partial T_{2}$ 
such that for each component $t$ of $T_1$ there exists 
a component $u$ of $T_2$ with $\partial t=\partial u$. 
 Two local moves $(T_{1},T_{2})$ and $(U_{1},U_{2})$
are {\it equivalent}, 
denoted by $(T_{1},T_{2})\cong (U_{1},U_{2})$, 
if there is an orientation preserving 
self-homeomorphism $\psi :B^{3}\rightarrow B^{3}$ such that $\psi (T_{i})$ 
and $U_{i}$ are ambient isotopic in $B^3$ relative to $\partial
B^{3}$ for 
$i=1,2$. 

A local move $(T_1,T_2)$ is {\em trivial}, if $(T_1,T_2)$ is 
equivalent to a local move $(T_1,T_1)$. 
Let $(T_1,T_2)$ be a local move, and let $t_1,t_2,...,t_k$ and 
$u_1,u_2,...,u_k$ be the components with 
$\partial t_i=\partial u_i\ (i=1,2,...,k)$ of 
$T_1$ and $T_2$ respectively.  
We call $(T_1,T_2)$ a {\em $k$-component Brunnian
local move} $(k\geq 2$), or {\em $B_k$-move}, if 
each local move $(T_1-t_i,T_2-u_i)$ is trivial $(i=1,2,...,k)$ 
\cite{T-Y0}. 
If $(T_{1},T_{2})$ is Brunnian, then $(T_{2},T_{1})$ is also Brunnian. 

Let $K_1$ and $K_2$ be oriented knots in the oriented three-sphere $S^3$. 
We say that $K_2$ is {\em obtained from $K_1$ by 
a local move $(T_{1},T_{2})$} if 
there is an orientation preserving embedding 
$h:B^{3}\rightarrow S^{3}$ such that $(h^{-1}(K_1),h^{-1}(K_2))
\cong (T_1,T_2)$ and $K_{1}-h(B^{3})=K_{2}-h(B^{3})$ together 
with orientations. Two oriented knots $K_{1}$ and
$K_{2}$ are {\it $B_{k}$-equivalent} if 
$K_{2}$ is obtained from $K_{1}$ by a finite sequence of $B_{k}$-moves 
and ambient isotopies. This relation is an equivalence relation on knots. 

We have the following geometric condition for knots to have 
same Vassiliev invariant. 

\medskip
{\bf Theorem 1.} (cf. Goussarov-Habiro Theorem 
\cite{Habiro1}, \cite{Gus2}) 
{\it Two knots $K_1$ and $K_2$ are $B_{l+1}$-equivalent if and only if 
their values of any Vassiliev invariant of order $\leq l-1$  
are equal. }

\medskip
{\bf Remark.} 
In \cite{Hab}, \cite{M-Y}, they independently show that 
the $B_{l+1}$-equivalence classes coincide with the $C_l$-equivalence classes. 
So the theorem above is same as Goussarov-Habiro Theorem.  
Although we can have Theorem 1 as a corollary of Goussarov-Habiro Theorem, 
we will give a self-contained proof of it. 
The author believes that it is worth to do so, because 
our proof is simple and short, and the known proofs 
of Goussarov-Habiro Theorem \cite{Habiro1}, \cite{Gus2}, \cite{T-Y2} 
are rather complicated.

%This does not mean that we give a short proof for Goussarov-Habiro Theorem. 
%In fact, both the proofs in \cite{Hab} and \cite{M-Y} 
%of that the $B_{l+1}$-equivalence implies 
%$C_l$-equivalence require some preliminaries lemmas shown in 
%\cite{Habiro1} and \cite{T-Y2} respectively. 
%And hence, in order to prove Goussarov-Habiro Theorem via Theorem 3, we 
%need such lemmas.

\medskip
Let $l$ be a positive integer and $k_1,...,k_l(\geq 2)$ integers. 
Suppose
that for each
$P\subset\{1,...,l\}$ an oriented knot $K_P$ in $S^3$ is assigned. Suppose 
that there are orientation preserving embeddings 
$h_i:B^3\rightarrow S^3$ $(i=1,...,l)$ such that \\
(1) $h_i(B^3)\cap h_j(B^3)=\emptyset$ if $i\neq j$,\\
(2) $K_P-\bigcup_{i=1}^l h_i(B^3)=K_{P'}-\bigcup_{i=1}^l h_i(B^3)$ together
with orientation 
for any subsets $P,P'\subset \{1,...,l\}$,\\
(3) $(h_i^{-1}(K_\emptyset),h_i^{-1}(K_{\{1,...,l\}}))$ is 
a $B_{k_i}$-move $(i=1,...,l)$, and\\
(4) $K_P\cap h_i(B^3)=\left\{
\begin{array}{ll}
K_{\{1,...,l\}}\cap h_i(B^3) & \mbox{if $i\in P$},\\
K_\emptyset\cap h_i(B^3) & \mbox{otherwise}.
\end{array}
\right.$\\
Then we call the set of oriented knots $\{K_P|P\subset\{1,...,l\}\}$ a 
{\em singular knot of type $B(k_1,...,k_l)$}. 
Let ${\cal K}$ be the set of knots, $A$ an abelian group, 
and $\varphi: {\cal K}\rightarrow A$ an invariant. 
We say that $\varphi$ is a {\em finite type invariant 
of type $B(k_1,...,k_l)$} if for a singular knot
$\{K_P|P\subset\{1,...,l\}\}$ of type 
$B(k_1,...,k_l)$,  
\[\sum_{P\subset\{1,...,l\}}(-1)^{| P|}\varphi(K_P)=0\in A.\]
Since a $B_2$-move is realized by some crossing changes 
we have that an invariant $\varphi:{\cal K}\rightarrow A$ is 
a finite type invariant of type $B(\underbrace{2,...,2}_{l})$ if and 
only if it is a Vassiliev invariant of order $\leq l-1$. 

In order to prove Theorem 1, we need the following theorems.

\medskip
{\bf Theorem 2.} (cf. \cite[Theorem 5.4]{Habiro1}) {\it The
$B_k$-equivalence classes, denoted by ${\cal K}/B_k$, of oriented knots in
$S^3$ forms an abelian group under connected sum of oriented knots.}

\medskip
{\bf Theorem 3.} (cf. \cite[Theorem 1.2]{T-Y2}) 
{\it Let $l(\geq2)$ and $k_1,...,k_l(\geq 2)$ be integers, and 
$k-1=(k_1-1)+\cdots+(k_l-1)$. Then the projection 
$p_k:{\cal K}\rightarrow{\cal K}/B_k$ is a finite type invariant 
of type $B(k_1,...,k_l)$.}

\medskip
{\bf Remark.} 
Since a $C_k$-move is a $B_{k+1}$-move,     
Theorem 2 follows \cite[Theorem 5.4]{Habiro1}. 
Theorem 3 is similar to \cite[Theorem 1.2]{T-Y2}. 
In order to give a self-contained proof of Theorem 1, we will give 
self-contained proofs of Theorems 2 and 3. 
Although the outlines of the proofs of Theorems 2 and 3 
are same as those in \cite{Habiro1} (or \cite{T-Y}, \cite{T-Y2}) 
and \cite{T-Y2} respectively, our proofs are simpler than theirs.

\bigskip\noindent
{\bf 2. Band description}

\bigskip
Let $(T_1,T_2)$ be $k$-component Brunnian local move. 
Let $T\subset B^3$ be the trivial $k$-string tangle illustrated in 
Fig 1, and $D$ a disjoint union of $k$ disks bounded by 
$T$ and arcs in $\partial B^3$ (see Fig 2). 
Note that $(T_1,T_2)$ is equivalent to a local move $(S,T)$. 
Then the pair $(S,\partial D-T)$ is said to be {\em $B_k$-link model} 
(see Fig 3). 

\begin{center} 
\begin{tabular}{cc}
\includegraphics[trim=0mm 0mm 0mm 0mm, width=.2\linewidth]
{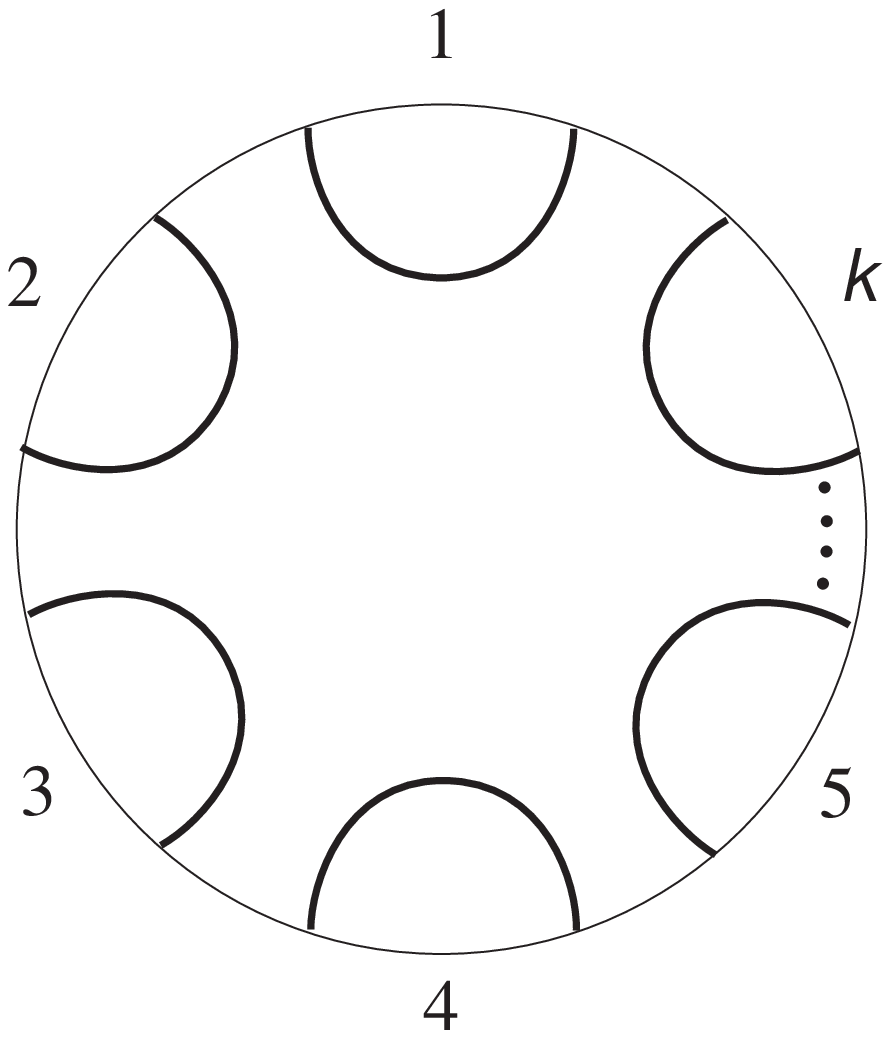} \hspace*{1in} &
\includegraphics[trim=0mm 0mm 0mm 0mm, width=.2\linewidth]
{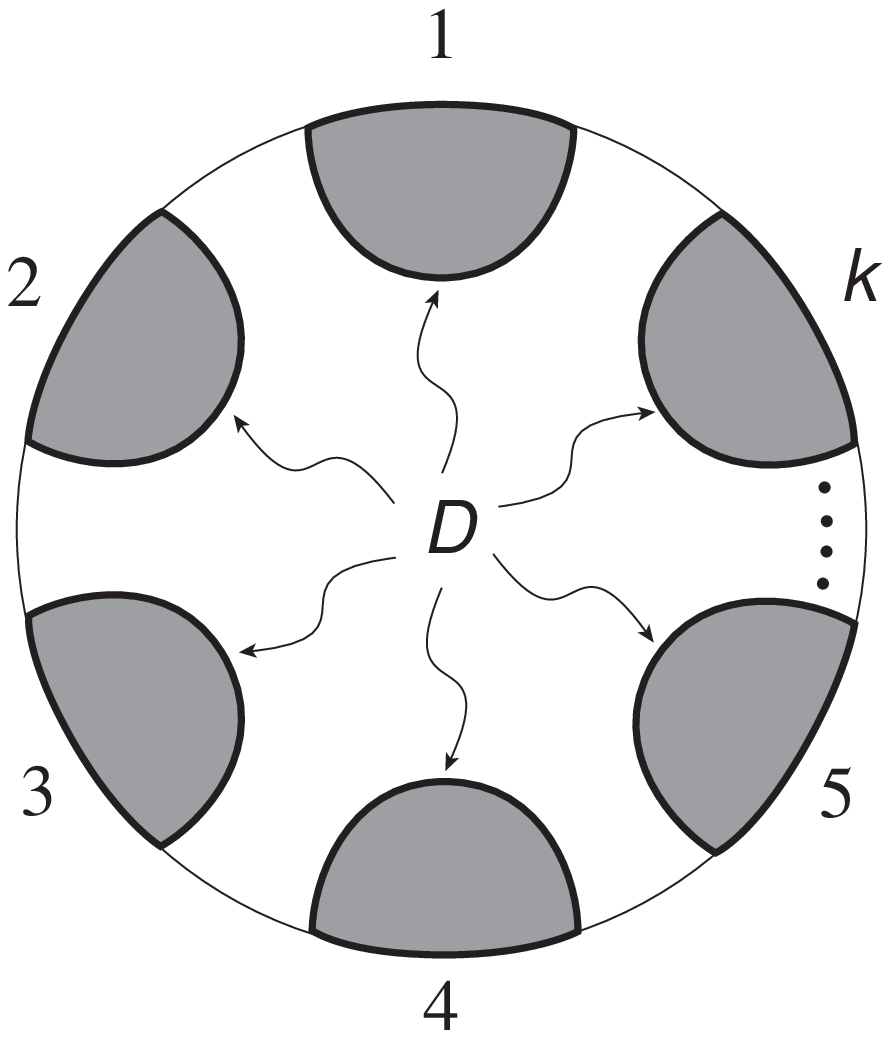}\\
Fig. 1 \hspace*{1in} &
Fig. 2
\end{tabular}
\end{center}

\medskip
\begin{center}
\includegraphics[trim=0mm 0mm 0mm 0mm, width=.7\linewidth]
{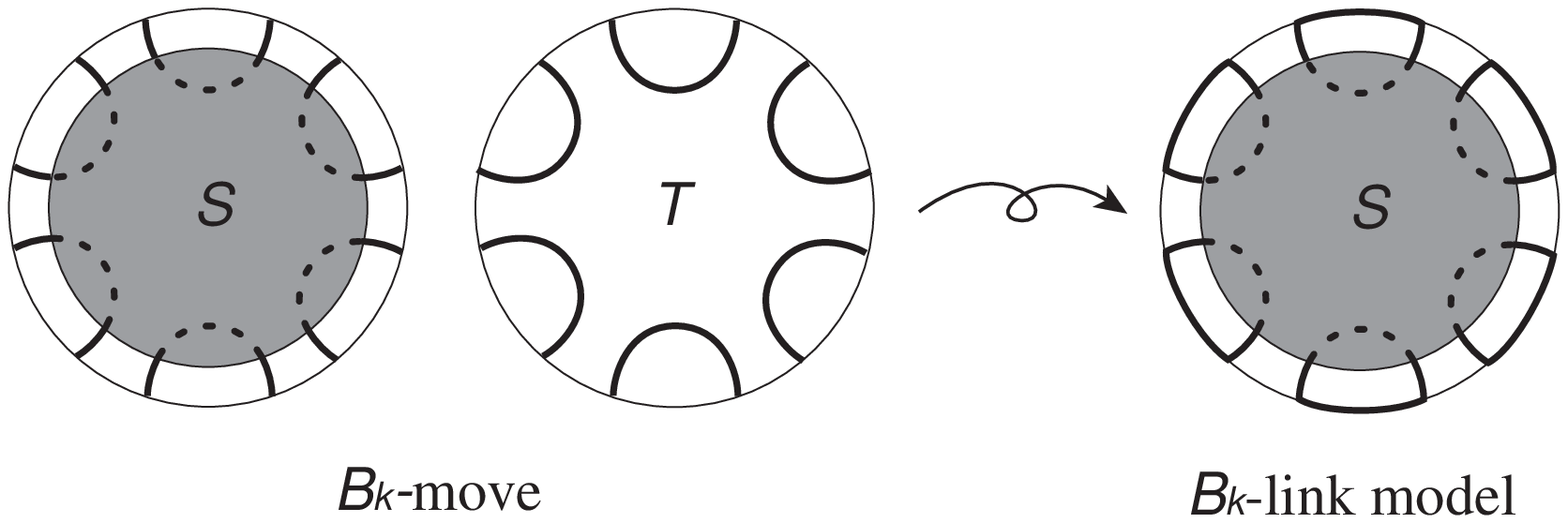}

Fig. 3
\end{center}
 
Let $(\alpha_{i},\beta_{i})$ be $B_{\rho(i)}$-link 
models $(i=1,...,l)$, and $K$ an oriented knot (resp. a tangle). Let
$\psi_{i}:B^{3}\rightarrow
S^3$ (resp. $\psi_{i}:B^{3}\rightarrow
{\rm int}B^3$)  be an orientation preserving embedding for $i=1,...,
l$ and 
$b_{1,1},b_{1,2},...,b_{1,\rho(1)},b_{2,1},b_{2,2},...,b_{2,\rho(2)},
...,b_{l,1},b_{l,2},...,b_{l,\rho(l)}$
mutually disjoint disks embedded in $S^3$ (resp. $B^3$). Suppose that 
they satisfy the following
conditions;\\
(1) $\psi_{i}(B^{3})\cap \psi_{j}(B^{3})=\emptyset$ if $i\neq j$,\\
(2) $\psi_{i}(B^{3})\cap K=\emptyset$  for each $i$,\\
(3) $b_{i,k}\cap K=\partial b_{i,k}\cap K$ is an arc for each $i,k$,\\
(4) $b_{i,k}\cap (\bigcup_{j=1}^{l} \psi_{j}(B^{3}))=
\partial b_{i,k}\cap \psi_{i}(B^{3})$ is a component of 
$\psi_{i}(\beta_{i})$ for each $i,k$.\\
Let $J$ be an oriented knot (resp. a tangle) defined by
\[
J=K\cup (\bigcup_{i,k}\partial b_{i,k})\cup 
(\bigcup_{i=1}^{l}\psi_{i}(\alpha_{i})) - 
\bigcup_{i,k}{\rm int}(\partial b_{i,k}\cap K) - 
\bigcup_{i=1}^{l}\psi_{i}({\rm int}\beta_{i}),
\]
where the orientation of $J$ 
coincides that of $K$ on $K-\bigcup_{i,k}b_{i,k}$ 
if $K$ is oriented.  We call each $b_{i,k}$ a {\it band}. 
Each image $\psi_{i}(B^{3})$ is called a {\it link ball}.  We set
${\cal B}_i=((\alpha_i,\beta_i),\psi_i,\{b_{i,1},...,b_{i,\rho(i)}\})$
and call ${\cal B}_i$ a {\it $B_{\rho(i)}$-chord}. We denote $J$  by
$J=\Omega(K;\{{\cal B}_1,...,{\cal B}_l\})$, and say that
$J$ is a {\em band sum} of
$K$ and chords ${\cal B}_1,...,{\cal B}_l$ or 
a {\em band sum} of
$K$ and $\{{\cal B}_1,...,{\cal B}_l\}$. 

\medskip
From now on we consider knots up to ambient isotopy of 
$S^3$ and tangles up to ambient isotopy of
$B^3$ relative to $\partial B^3$ without explicit mention. 

By the definitions of $B_k$-move and $B_k$-link model,  
we have 

\medskip
{\bf Sublemma 4.} (cf. \cite[Sublemmas 3.3 and 3.5]{T-Y2}) 
{\it {\rm (1)} A local move $(T_1,T_2)$ is a $B_k$-move if and only if 
$T_1$ is a band sum of $T_2$ and a $B_k$-link model. 

{\rm (2)} A knot $J$ is obtained from a knot $K$ by a single $B_k$-move 
if and only if $K$ is a band sum of $J$ and a $B_k$-link model. 
$\Box$}

\medskip
Note that, by Sublemma 4 (1), a set of knots {\bf K} is 
a singular knot of type $B(k_1,...,k_l)$ if and only if 
there is a knot $K$ and a band sum 
$J=\Omega(K;\{{\cal B}_1,...,{\cal B}_{l}\})$ 
 of $K$ and $B_{k_i}$-chords ${\cal B}_i$ $(i=1,...,l)$ such that 
\[{\mathrm\bf K}=
\left\{\left. \Omega\left(K;\bigcup_{i\in P}\{{\cal B}_{i}\}\right) \right| 
P\subset\{1,...,l\}\right\}. \] 

\medskip
{\bf Sublemma 5.} (cf. \cite[Sublemma 3.5]{T-Y2}) 
 {\it Let $K$, $J$ and $I$ be oriented
knots $($or tangles$)$. Suppose that
$J=\Omega(K;\{{\cal B}_1,...,{\cal B}_l\})$ for some chords ${\cal
B}_1,...,{\cal B}_l$
and $I=\Omega(J;\{{\cal B}\})$ for some $B_k$-chord ${\cal B}$. Then there
is a $B_k$-chord ${\cal B}'$ such that 
$I=\Omega(K;\{{\cal B}_1,...,{\cal B}_l,{\cal B}'\})$. 
Moreover, if a subset $P$ of $\{1,...,l\}$ satisfies that 
the link ball or the bands of $\cal B$ intersect either the link ball 
or the bands of ${\cal B}_j$ only if $i\in P$, 
then $\Omega(\Omega(K;\bigcup_{i\in P}\{{\cal B}_i\});\{{\cal B}\})=
\Omega(K;(\bigcup_{i\in P}\{{\cal B}_i\})\cup\{{\cal B}'\})$.}

\medskip
{\bf Proof.} If the bands and the
link ball of ${\cal B}$ are disjoint from those of ${\cal B}_1,...,{\cal
B}_l$ then we have that 
$I=\Omega(K;\{{\cal B}_1,...,{\cal B}_l,{\cal B}\})$. If not then we
deform $I$ up to
ambient isotopy as follows. 
By thinning and shrinking the bands and the link
ball of ${\cal B}$ respectively, we may assume that the link
ball of ${\cal B}$ intersects neither 
the bands nor the link balls of ${\cal B}_1,...,{\cal B}_l$. 
And by sliding the bands of ${\cal B}$ along $J$, we may also assume that 
the intersection of the bands with $J$ is 
disjoint from the bands and the link balls of ${\cal B}_1,...,{\cal B}_l$. 
Then we
sweep the bands of
${\cal B}$ out of the link balls of ${\cal B}_1,...,{\cal B}_l$. Note
that this is
always possible since the tangles of a local move are trivial. 
Finally we sweep the
intersection of the bands of
${\cal B}$ and the bands of ${\cal B}_1,...,{\cal B}_l$ out of the
intersection of the bands
of ${\cal B}_1,...,{\cal B}_l$ and $K$. Let ${\cal B}'$ be the result
of the deformation of
${\cal B}$ described above. Then it is not hard to see that 
${\cal B}'$ is a desired chord. $\Box$

\medskip
By repeated applications of Sublemmas 4 and 5 we immediately 
have the following lemma.

\medskip
{\bf Lemma 6.} (cf. \cite[Lemma 3.6]{T-Y2}) 
 {\it Let $k$ be a positive integer 
and let $K$ and $J$ be oriented knots $($or tangles$)$. 
Then $K$ and $J$ are $B_{k}$-equivalent if and only if $J$ is a band
sum of $K$ and some $B_{k}$-link models. $\Box$}

\medskip 
Since the local moves illustrated in Figs 4 and 5 
are $B_{k+1}$-move and $B_{j+k-1}$-move respectively, 
the following two lemmas follow Sublemma 5. 

\begin{center} 
\begin{tabular}{cc}
\hspace*{-10mm} \includegraphics[trim=0mm 0mm 0mm 0mm, width=.35\linewidth]
{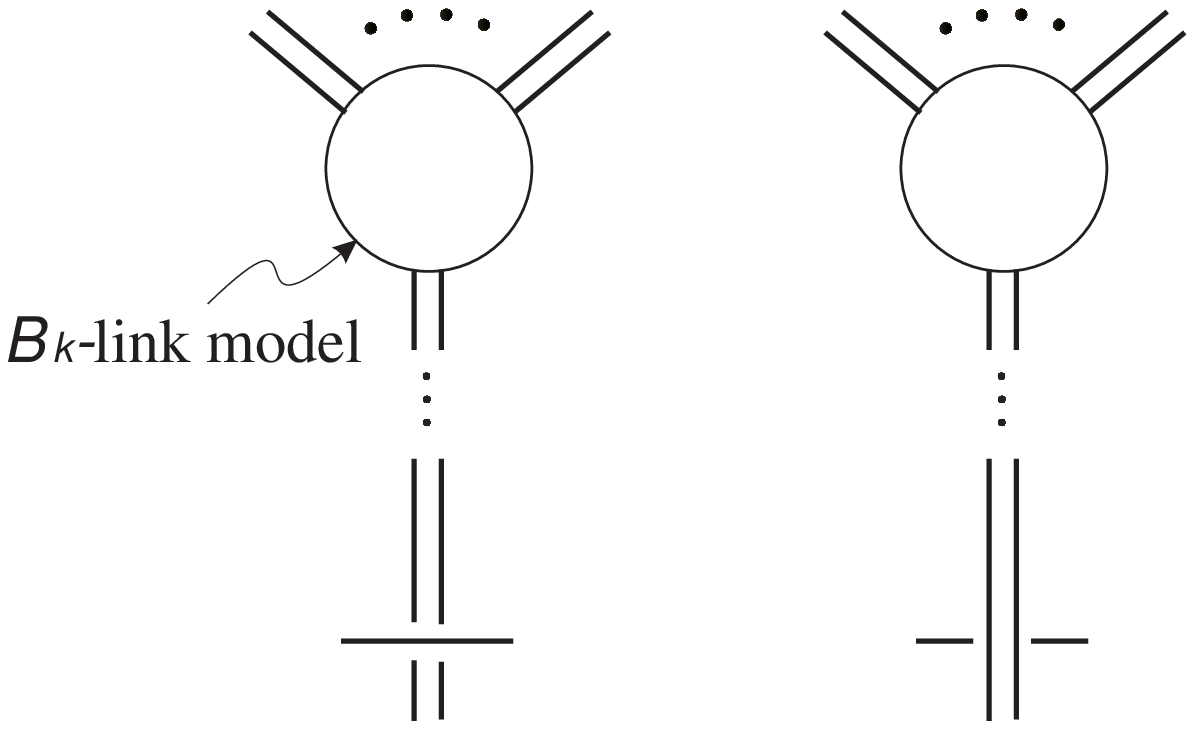}  &
\includegraphics[trim=0mm 0mm 0mm 0mm, width=.6\linewidth]
{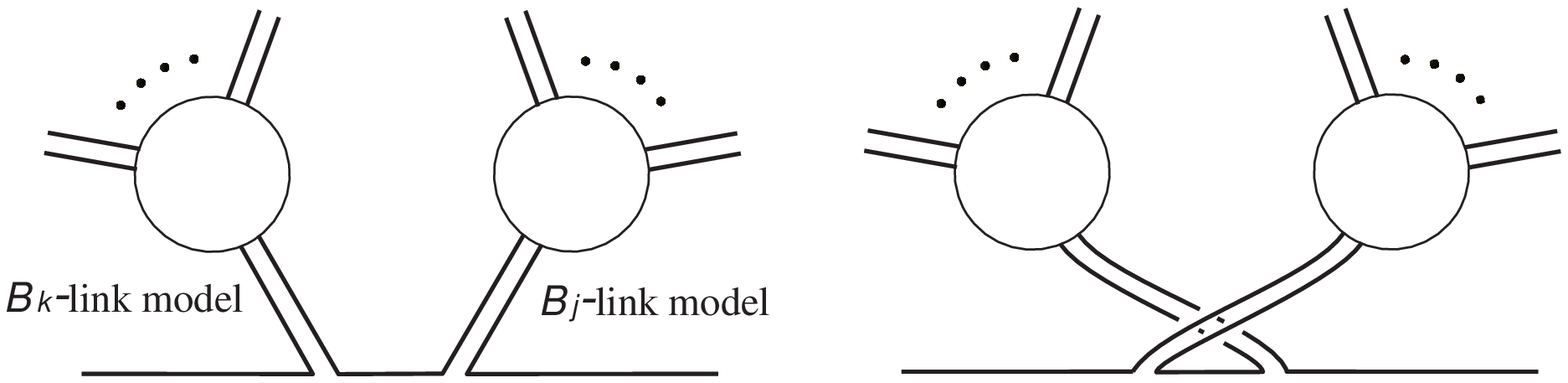}\\
Fig. 4 &
Fig. 5
\end{tabular}
\end{center}

\medskip
{\bf Lemma 7.} (cf. \cite[Lemma 3.8]{T-Y2}) 
{\it  Let $K$, $J=\Omega(K;\{{\cal B}_1,...,{\cal B}_l,{\cal B}_0\})$ 
 and $I=\Omega(K;\{{\cal B}_1,...,{\cal B}_l,{\cal B}'_0\})$ be 
 oriented knots, where 
${\cal B}_1,...,{\cal B}_l$ are chords and ${\cal B}_0,{\cal B}'_0$ 
are $B_k$-chords. Suppose that $J$ and $I$ differ locally as illustrated 
in Fig {\rm 4}, i.e., 
$I$ is obtained from $J$ by a crossing change 
between $K$ and a band of ${\cal B}_0$. Then $I$ is obtained from $J$ 
by a $B_{k+1}$-move. Moreover, there is a $B_{k+1}$-chord $\cal B$ such 
that $\Omega(K;(\bigcup_{i\in P}\{{\cal B}_i\})\cup\{{\cal B}_0\})=
\Omega(K;(\bigcup_{i\in P}\{{\cal B}_i\})\cup\{{\cal B}'_0,{\cal B}\})$ 
for any subset $P$ of $\{1,...,l\}$. $\Box$}

\medskip
{\bf Lemma 8.} (cf. \cite[Lemma 3.9]{T-Y2}) 
{\it  Let $K$, 
$J=\Omega(K;\{{\cal B}_1,...,{\cal B}_l,{\cal B}_{0j}, 
{\cal B}_{0k}\})$ and $I=\Omega(K;\{{\cal B}_1,...,{\cal B}_l,
{\cal B}'_{0j},{\cal B}'_{0k}\})$ be oriented knots, where 
${\cal B}_1,...,{\cal B}_l$ are chords and 
${\cal B}_{0j},{\cal B}'_{0j}$ $($resp. ${\cal B}_{0k},{\cal B}'_{0k})$ 
are $B_j$-chords $($resp. $B_k$-chords$)$. 
Suppose that $J$ and $I$ differ locally as illustrated in Fig {\rm 5}. 
Then $I$ is obtained from $J$ by a $B_{j+k-1}$-move. 
Moreover, there is a $B_{j+k-1}$-chord $\cal B$ such 
that $\Omega(K;(\bigcup_{i\in P}\{{\cal B}_i\})\cup\{{\cal B}_{0j}, 
{\cal B}_{0k}\})=\Omega(K;(\bigcup_{i\in P}\{{\cal B}_i\})\cup
\{{\cal B}'_{0j},{\cal B}'_{0k},{\cal B}\})$ 
for any subset $P$ of $\{1,...,l\}$. $\Box$}

\medskip\noindent
We call the change from $J$ to $I$ in Lemma 8 a {\em band exchange}.

For $C_k$-move, \lq band description' is also defined, and 
Sublemmas 4, 5, Lemmas 6, 7 and 8 hold \cite{T-Y2}. 
But the proofs are not obvious as ours. 
In fact, little complicated arguments are needed. 
In contrast, we need some arguments to prove the following lemma, which 
is trivial for $C_k$-move. 

\medskip
{\bf Lemma 9.} {\it  Let $(T_1,T_2)$ be a $B_k$-move. 
For any integer $l(\leq k)$, $T_2$ is obtained from $T_1$ by $B_l$-moves. 
In particular, $B_k$-equivalence knots are $B_l$-equivalent.} 

\medskip
{\bf Proof.} 
Let $t_1,...,t_k$ and $u_1,...,u_k$ be the components 
with $\partial t_i=\partial u_i\ (i=1,...,k)$ of
$T_1$ and $T_2$ respectively. 
We may assume that $(T_1,T_2)$ has a diagram in the unit disk such that 
both $T_1-t_1$ and $T_2$ has no crossings. 

Since $(T_1-t_2,T_2-u_2)$ is a trivial local move, 
$T_2$ is obtained from $T_1$ by $B_2$-moves that correspond to crossing 
changes between $t_1$ and $t_2$. By Lemma 6, $T_1$ is a band sum 
$\Omega(T_2;{\bf B}_2)$ of $T_2$ and a set ${\bf B}_2$ of 
$B_2$-chords. Note that any bands of the $B_2$-chords does not 
intersect to $T_2-(u_1\cup u_2)$. 

Since $(\Omega(T_2;{\bf B}_2)-t_3,T_2-u_3)=(T_1-t_3,T_2-u_3)$ 
is a trivial local move, 
$T_2$ is obtained from $T_1$ by $B_3$-moves that correspond to crossing 
changes between $t_3$ and some bands of the $B_2$-chords. 
By Lemma 6, $T_1$ is a band sum 
$\Omega(T_2;{\bf B}_3)$ of $T_2$ and a set ${\bf B}_3$ of 
$B_3$-chords. Note that any bands of the $B_3$-chords does not 
intersect to $T_2-(u_1\cup u_2\cup u_3)$. 

Repeat these processes, then we have the conclusion. $\Box$

\bigskip\noindent
{\bf 3. Proofs of Theorems 1, 2 and 3.}

\bigskip
\noindent
{\bf Proof of Theorem 3.} 
Let $k_1,...,k_l(\geq 2)$ be positive integer and 
$k-1=(k_1-1)+\cdots+ (k_l-1)$. 
Let $K_0$ be a knot and $K_1$ a band sum of $K_0$ and $B_{k_j}$-chords 
${\cal B}_{k_j,j}$ $(j=1,...,l)$. 
It is sufficient to show that 
\[\sum_{P\subset\{1,...,l\}}(-1)^{|
P|} \left[\Omega\left(K_0;\bigcup_{j\in P}\{{\cal B}_{k_j,j}\}\right)
\right]
=0\in {\cal K}/B_k,\]
where $[K]$ is the $B_k$-equivalent class which contains a knot $K$. 

Set
\[K_P=\Omega\left(K_0;\bigcup_{j\in P}\{{\cal
B}_{k_j,j}\}\right).\]

\medskip
{\bf Claim.} The knot $K_1(=K_{\{1,...,l\}})$ is $B_k$-equivalent to a band sum 
%\[\Omega\left(K_0;\bigcup_{i,j}\{{\cal B}_{i,j}\}\right)\]
of $K_0(=K_{\emptyset})$ and a set of {\em local chords} $\bigcup_{i,j}\{{\cal B}_{i,j}\}$
such that (1) ${\cal B}_{i,j}$ is a $B_i$-chord $(i<k)$ 
and it has an associated subset $\omega({\cal B}_{i,j})\subset\{1,...,l\}$ 
with $\sum_{t\in\omega({\cal B}_{i,j})}(k_t-1)\leq i-1$, and 
(2) for each $P\subset\{1,...,l\}$
\[K_P=\Omega\left(K_0;\bigcup_{\omega({\cal B}_{i,j})
\subset P}\{{\cal B}_{i,j}\}\right).\]
Here a chord ${\cal B}_{i,j}$ is called a {\it local chord} 
if there is a 3-ball $B$ such that $B$ contains all 
of the bands and the link ball of 
${\cal B}_{i,j}$, $B$ does not intersect any other bands and link
balls, and that $(B,B\cap K_0)$ is a trivial ball-arc pair. 

\medskip
Before proving Claim, we will finish the proof of Theorem 3. 
Suppose $K_1$ is $B_k$-equivalent to a band sum of $K_0$
and some local chords ${\cal B}_{ij}$'s. 
Such a local chord ${\cal B}_{ij}$ represents 
a knot $K_{ij}$ connected summed to $K_0$. 
So the band sum is a connected sum of $K_0$ and $K_{ij}$'s. 
Then we have 
\[\begin{array}{rl}
 &\displaystyle{
  \sum_{P\subset\{1,...,l\}}
 (-1)^{| P|}\left[
 \Omega\left(K_0;\bigcup_{\omega({\cal B}_{i,j})\subset P}
 \{{\cal B}_{i,j}\}\right)\right]}\\
 =&\displaystyle{
  \sum_{P\subset\{1,...,l\}}(-1)^{| P|}
  \left([K_0]
 +\sum_{\omega({\cal B}_{i,j})\subset P}[K_{i,j}]\right)}\\
=&\displaystyle{
 \sum_{P\subset\{1,...,l\}}(-1)^{| P|}[K_0]+
 \sum_{P\subset\{1,...,l\}}(-1)^{| P|}\left(\sum_{\omega
 ({\cal B}_{i,j})\subset P}[K_{i,j}]\right)}\\
=&\displaystyle{
 0+\sum_{i,j}\left(\sum_{P\subset\{1,...,l\},\omega({\cal B}_{i,j}
) \subset P}(-1)^{| P|}\right)[K_{i,j}]}.
\end{array}\]
We consider the coefficient of $[K_{i,j}]$. 
Since $\sum_{t\in\omega({\cal B}_{i,j})}(k_t-1)< k-1$, 
$\omega({\cal B}_{i,j})$ is a proper subset of
$\{1,...,l\}$. We may assume that $\omega({\cal B}_{i,j})$ 
does not contain $a\in\{1,...,l\}$.
Then we have 
\[\begin{array}{rl}
\displaystyle{ 
\sum_{P\subset\{1,...,l\},\omega({\cal B}_{i,j})\subset P}
(-1)^{| P|}}=
&\displaystyle{
\sum_{P\subset\{1,...,l\}\setminus\{a\},\omega({\cal B}_{i,j})
\subset P}(-1)^{| P|}}\\
&\hspace*{3em}+
\displaystyle{\sum_{P\subset\{1,...,l\}\setminus\{a\},\omega({\cal
B}_{i,j})\subset P}(-1)^{| P\cup\{a\}|}=0.}
\end{array}\]
Thus, we have the conclusion.  $\Box$

\medskip
Now we will show Claim. 

\medskip
{\bf Proof of Claim.} 
We first set $\omega({\cal B}_{k_j,j})=\{j\}$ for $j=1,...,l$. 
Then we have $\sum_{t\in\omega({\cal B}_{k_j,j})}(k_t-1)=k_j-1<k-1$ and 
\[K_P=\Omega\left(K_0;\bigcup_{\omega({\cal
B}_{k_j,j})\subset P}\{{\cal B}_{k_j,j}\}\right).\] 
Note that a crossing change between bands 
can be realized by crossing changes
between $K_0$ and a band as illustrated in Fig. 6. Therefore we can
deform each chord into a local
chord by (i) crossing changes between $K_0$ and bands, and 
(ii) band exchanges.

(i)~~When we perform a crossing change between $K_0$ and a $B_p$-band of a
$B_p$-chord ${\cal B}_{p,q}$ with $p\leq k-2$, by using Lemma 7, 
we introduce a new 
$B_{p+1}$-chord ${\cal B}_{p+1,r}$ and we set $\omega({\cal
B}_{p+1,r})=\omega({\cal B}_{p,q})$ so that the conditions (1) and (2) 
still holds. 
By Lemma 7, a crossing change between $K_0$ and a
$B_{k-1}$-band is realized by a 
$B_{k}$-move and therefore does not change the
$B_{k}$-equivalence class. 

(ii)~~When we perform a band 
exchange between a $B_p$-chord ${\cal B}_{p,q}$ and a $B_r$-chord ${\cal
B}_{r,s}$ with $p+r\leq k$, then, by using Lemma 8,  
we introduce a new
$B_{p+r-1}$-chord ${\cal B}_{p+r-1,n}$ and set $\omega({\cal
B}_{p+r-1,n})=\omega({\cal B}_{p,q})\cup\omega({\cal B}_{r,s})$ so that the
conditions (1) and (2) still holds. 
By Lemmas 8 and 9, 
a band exchange between a $B_p$-chord ${\cal B}_{p,q}$ and a $B_r$-chord ${\cal
B}_{r,s}$ with $p+r\geq k+1$ does not change the $B_{k}$-equivalence class. $\Box$

\begin{center} 
\includegraphics[trim=0mm 0mm 0mm 0mm, width=.7\linewidth]
{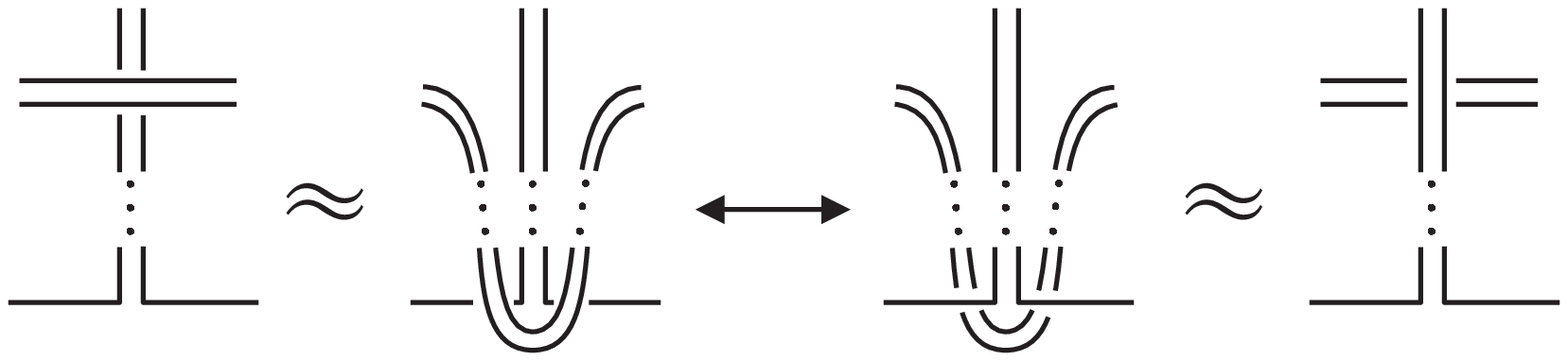}\\
Fig. 6
\end{center}

\medskip
{\bf Proof of Theorem 2.} 
It is sufficient to show that the existence of the inverse
element. Let $K$ be a knot. Suppose that there is a 
knot $J$ such that $K\#J$ is $B_{k}$-equivalent to a trivial 
knot $O$. Then, by Lemma 6, we have that
$O$ is a band sum of $K\#J$ and some $B_{k}$-chords. 
By using Lemma 7, we deform $O$ up to
$B_{k+1}$-equivalence so that the $B_k$-chords are local chords. 
Then we have that the result is a connected sum of $K\#J$ 
and some knots $K_1,...,K_n$ that correspond to the local chords. 
Namely $K\#J\#K_1\#\cdots\# K_n$ is $B_{k+1}$-equivalent to $O$. 
Thus $J\#K_1\#\cdots\# K_n$ is the desired knot. $\Box$

\medskip
{\bf Proof of Theorem 1.} 
It is not hard to see that a $B_{l+1}$-equivalence knots are  
$l$-similar \cite{Tan} ($(l-1)$-equivalent \cite{Gus1}). 

By Theorem 3, the projection 
$p_{l+1}:{\cal K}\rightarrow {\cal K}/B_{l+1}$ is a Vassiliev 
invariant of order $\leq l-1$. 
If two knots have same values of any 
Vassiliev invariant of order $\leq l-1$, then they are 
$B_{l+1}$-equivalent.  $\Box$

\bigskip
\footnotesize{
 }

\end{document}